\documentclass[twocolumn]{autart}    % Enable this line and disable the
\usepackage{graphicx}          % Include this line if your
\usepackage{amsbsy}
\usepackage{booktabs}
\usepackage{amssymb}
\usepackage{subeqnarray}
\begin{document}

\begin{frontmatter}
%\runtitle{Insert a suggested running title}  % Running title for regular
                                              % papers but only if the title
                                              % is over 5 words. Running title
                                              % is not shown in output.

\title{Sequential Fusion Estimation for Clustered Sensor Networks} % Title, preferably not more  \thanksref{footnoteinfo}
                                                % than 10 words.

%\thanks[footnoteinfo]{This paper was not presented at any IFAC
%meeting. Corresponding author W.~A.~Zhang. Tel. +86-571-85290531.}

\author[Paestum]{Wen-An Zhang}\ead{wazhang@zjut.edu.cn},    % Add the
\author[Rome]{Ling Shi}\ead{eesling@ust.hk},               % e-mail address
%\author[Baiae]{Publius Maro Vergilius}\ead{vergilius@culture.ir}  % (ead) as shown

\address[Paestum]{Department of Automation, Zhejiang University of Technology, Hangzhou 310023, China}  % Please supply
\address[Rome]{Electronic and Computer Engineering, Hong Kong University of Science and Technology, Clear Water Bay, Kowloon, Hong Kong}             % full addresses
%\address[Baiae]{The White House, Baiae}        % here.

                             % keyword list or with the
                                          % help of the Automatica
                                          % keyword wizard

\begin{abstract}                          % Abstract of not more than 200 words.
We consider multi-sensor fusion estimation for clustered sensor networks. Both sequential measurement fusion and state fusion estimation methods are presented. It is shown that the proposed sequential fusion estimation methods achieve the same performance as the batch fusion one, but are more convenient to deal with asynchronous or delayed data since they are able to handle the data that is available sequentially. Moreover, the sequential measurement fusion method has lower computational complexity than the conventional sequential Kalman estimation and the measurement augmentation methods, while the sequential state fusion method is shown to have lower computational complexity than the batch state fusion one. Simulations of a target tracking system are presented to demonstrate the effectiveness of the proposed results.
\end{abstract}

\end{frontmatter}

\section{Introduction}
Fusion estimation for sensor networks has attracted much
research interest during the last decade, and has found applications
in a variety of areas \cite{IlicMD10,CaoXH14,ChenJM13,HeX11,OkaA10}. Compared with the
centralized structure, the distributed structure is more preferable for sensor
networks because of its reliability, robustness and low
requirement on network bandwidth \cite{HeX11,DongHL12,MillanP13}. When the number
of sensors is large, it is wasteful to embed in each sensor
an estimator and the communication burden is high.
Moreover, for long-distance deployed sensors, it may not be
possible to allocate communication channels for all sensors.
An improvement is to adopt the hierarchical structure for distributed estimation \cite{SongHY14},
by which all the sensors in the network are divided into several clusters
and the sensors within the same cluster are connected to a
cluster head (CH) which acts as a local estimator. Then, the distributed estimation is carried out in two stages.
In the first stage, the local estimator in each cluster fuses the measurements from its cluster to generate a local estimate. Then,
the local estimators exchange and fuse local estimates to produce a fused estimate to eliminate any disagreements among themselves.

\begin{figure}
\begin{center}
\includegraphics[height=3.0cm]{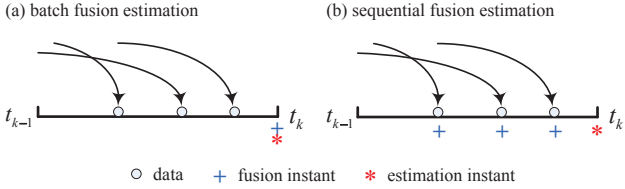}
\caption{Examples of batch fusion estimation and sequential fusion estimation.}
\label{fig1}
\end{center}
\end{figure}

Various results on multi-sensor fusion estimation for sensor networks have been available in the literature, including centralized fusion and distributed fusion, as well as measurement fusion and state fusion \cite{Roeckerja88,BarShalom95,Sunsl04,Deng06,Songeb07,Julier09,YYhu10,Zhang14,XiaYQ09,XingZR16}. However, most of the results are based on the idea of batch fusion, that is, measurements or local estimates are fused all at a time at the fusion instant until all of them are available at the estimator, as illustrated in Fig.\ref{fig1}(a). Such a batch fusion estimation may induce long computation delay, thus it is not appropriate for real-time applications. A possible improvement is to adopt the idea of sequential fusion, by which the measurements or local estimates are fused one by one according to the time order of the data arriving at the estimator, as illustrated in Fig.\ref{fig1}(b). In this way, the fusion and the state estimation could be carried out over the entire estimation interval, which help reduce computation burdens at the estimation instant and ultimately reduce the computation delay. Moreover, asynchronous or delayed data can be easily handled. Some relevant results on sequential fusion estimation have been presented in \cite{YanLP13} and \cite{Deng12}. The idea in \cite{YanLP13} is similar to the conventional sequential Kalman filtering, where the state estimate is updated several times by sequentially fusing the various measurements, and both the procedures of state prediction and measurement updating are involved over each step of the estimate updating, which incurs significant much computation cost. An alternative approach is to fuse all the measurements first, and then generate the state estimate based on the fused measurement. This is the novel method introduced in this paper. In \cite{Deng12}, the sequential covariance intersection (CI) fusion method was presented for state fusion estimation. However, the CI fusion is not optimal since the cross-covariances among the various local estimates are ignored.

In this paper, both sequential measurement fusion (SMF) estimation and state fusion estimation (SSF) methods are developed for clustered sensor networks, where the SMF is presented for local estimation, while the SSF is presented for state fusion estimation among all the local estimators. The main contributions of the paper are summarized as follows:\\
$~~$ 1) We present a design method for the SMF estimators. We show that the SMF estimator is equivalent to the conventional sequential Kalman (SK) and the batch measurement fusion (BMF) estimators, and is equivalent to the one designed based on measurement augmentation (MA). We also show that the SMF estimator has lower computational complexity than the estimators based on SK and MA.\\
$~~$ 2) We present a design method for the SSF estimators with matrix weights. We further show that the SSF estimator is equivalent to the batch state fusion (BSF) estimators with matrix weights but has much lower computational complexity.\\

\begin{figure}
\begin{center}
\includegraphics[height=4.8cm]{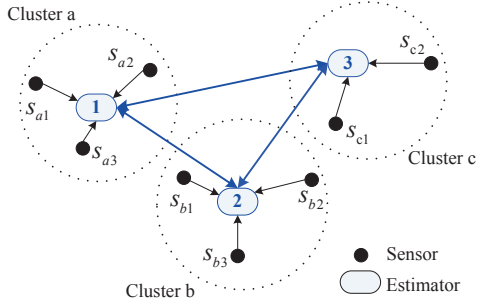}
\caption{A structure of hierarchical fusion estimation for clustered sensor networks.}
\label{fig2}
\end{center}
\end{figure}

\section{Problem Statement}
Consider the hierarchical fusion estimation for clustered sensor networks as shown in Fig.\ref{fig2},
where the plant, whose state is to be estimated, is described by the following discrete-time state-space model
\begin{eqnarray}
\label{eq:crosseq-1}
x(k+1)=A(k)x(k)+B(k)\omega(k)
\end{eqnarray}
where $x(k)\in\Re^{n_{x}}$ is the system state, and $\omega(k)\in\Re^{n_{\omega}}$ is a zero-mean white Gaussian noise with variance $Q_{\omega}$. A sensor network with $m$ clusters is deployed to monitor the state of system (\ref{eq:crosseq-1}). The set of the clusters is denoted by $\Phi=\{1,\ldots,m\}$. Let $\mathcal{N}_{s}=\{1,\ldots,n_{s}\}$ denote the $s$th cluster in the sensor network, where $s\in \Phi$ and $n_{s}$ is the number of sensors in the cluster $\mathcal{N}_{s}$. The $n_{s}$ sensors are connected to a cluster head (CH) $e_{s}$ serving as an estimator. The measurement equation of each sensor is given by
\begin{eqnarray}
\label{eq:crosseq-2}
y_{s,i}(k)=C(k)x(k)+\upsilon_{s,i}(k),~i\in\mathcal{N}_{s},~s\in \Phi
\end{eqnarray}
where $y_{s,i}(k)\in\Re^{q}$, $\upsilon_{s,i}(k)$ is a zero-mean white Gaussian noise with variance $R_{s,i}$, and $\upsilon_{s,i}(k)$ are mutually uncorrelated and are uncorrelated with $\omega(k)$.

As shown in Fig.\ref{fig2}, the fusion estimation is carried out in two stages. At the first stage, each CH collects and fuses measurements sequentially from its cluster, then generates a local estimate using the fused measurement. At the second stage, each CH collects local estimates from itself and the other CHs to produce a fused state estimate using the SSF method to improve estimation performance and eliminate any disagreements among the estimators.

\section{Design of the SMF Estimators}
This section is devoted to the design of the SMF estimators for each cluster. Consider cluster $\mathcal{N}_{s}$, $s\in \Phi$. For notational convenience, the subscript $s$ in the notations will be dropped in the remaining of this section, for example, $y_{s,i}$ is denoted as $y_{i}$ and $n_{s}$ is denoted as $n$. Denote $y_{f}^{s}$ as the fused measurement and $Y(k)=\{y_{1}(k),\ldots,y_{n}(k)\}$ as the set of measurements for fusion. Then it can be seen from Fig.\ref{fig1}(b) that $y_{f}^{s}$ is obtained by sequentially fusing the $n$ measurements. The fused measurement and its noise variance of the $j$th fusion over the interval $(k-1,k]$ is denoted by $y_{(j)}(k)$ and $R_{(j-1)}(k)$, respectively, where $j\in\{1,2,\ldots,n-1\}$. Denote the measurement noise of $y_{(j)}(k)$ as $\upsilon_{(j)}(k)$, then $y_{(j)}(k)=C(k)x(k)+\upsilon_{(j)}(k)$, and $R_{(j-1)}(k)=\mathrm{Cov}\{\upsilon_{(j)}(k)\}$. We now introduce the first main result on SMF estimator.

{\bf Theorem 1}. For the measurements in $Y(k)$, the SMF estimator is given by the following equations
\begin{eqnarray}
\label{eq:crosseq-3}
R_{(j)}(k)&=&\left[R_{(j-1)}^{-1}(k)+R_{j+1}^{-1}(k)\right]^{-1}\\
\label{eq:crosseq-4}
y_{(j)}(k)&=&R_{(j)}(k)\left[R_{(j-1)}^{-1}(k)y_{(j-1)}(k)\right.\nonumber\\
&&\left.+R_{j+1}^{-1}(k)y_{j+1}(k)\right]
\end{eqnarray}
where $j=1,\ldots,n-1$, $y_{(0)}(k)=y_{1}(k)$, $R_{(0)}(k)=R_{1}(k)$, and the fused measurement $y_{f}^{s}(k)$ and its noise variance $R_{f}^{s}(k)$ are given by $y_{f}^{s}(k)=y_{(n-1)}(k)$ and $R_{f}^{s}(k)=R_{(n-1)}(k)$, respectively. Moreover, one has $R_{(j)}(k)\leq R_{(j-1)}(k)$ and $R_{f}^{s}(k)\leq R_{i}(k)$, $i\in\{1,\ldots,n\}$.

{\bf Proof}. For brevity, the notation $k$ will be dropped in the following developments. Denote $f_{m}$ as the sequential measurement fusion operator, then $y_{(j)}=f_{m}\{y_{(j-1)},y_{j+1}\}$. Augment $y_{(j-1)}$ and $y_{j+1}$ to get
\begin{eqnarray}
\label{eq:crosseq-5}
z_{(j)}=\left[\begin{array}{c}y_{(j-1)}\\y_{j+1}\end{array}\right]=eCx+\bar{\upsilon}_{(j)}
\end{eqnarray}
where $e=[I~~I]^{\mathrm{T}}$ and $\bar{\upsilon}_{(j)}=[\upsilon_{(j-1)}^{\mathrm{T}}~~\upsilon_{j+1}^{\mathrm{T}}]^{\mathrm{T}}$. Let $\bar{R}_{(j)}=\mathrm{Cov}\{\bar{\upsilon}_{(j)}\}$. The term $z_{(j)}$ can be regarded as a measurement of $Cx$ with the measurement noise $\bar{\upsilon}_{(j)}$ and the measurement matrix $e$. Then by the weighted least square (WLS) estimation method, a least norm estimate of $Cx$ is given by
\begin{eqnarray}
\label{eq:crosseq-6}
\hat{z}_{(j)}=[e^{\mathrm{T}}\bar{R}_{(j)}^{-1}e]^{-1}e^{\mathrm{T}}\bar{R}_{(j)}^{-1}z_{(j)}
\end{eqnarray}
Since $\upsilon_{(j-1)}$ is related to $\{\upsilon_{1},\ldots,\upsilon_{j}\}$, it is uncorrelated with $\upsilon_{j+1}$. Thus, one has
\begin{eqnarray}
\label{eq:crosseq-7}
\bar{R}_{(j)}=\mathrm{diag}\{R_{(j-1)},R_{j+1}\}
\end{eqnarray}
Substituting (\ref{eq:crosseq-5}) and (\ref{eq:crosseq-7}) into (\ref{eq:crosseq-6}) yields
\begin{eqnarray}
\label{eq:crosseq-8}
\hat{z}_{(j)}=[R_{(j-1)}^{-1}+R_{j+1}^{-1}]^{-1}[R_{(j-1)}^{-1}y_{(j-1)}+R_{j+1}^{-1}y_{j+1}]
\end{eqnarray}
It can be seen from (\ref{eq:crosseq-8}) that $\hat{z}_{(j)}$ is a linear combination of $y_{(j-1)}$ and $y_{j+1}$, thus it can be regarded as the fused measurement $y_{(j)}$, that is, $y_{(j)}=\hat{z}_{(j)}$, which leads to equation (\ref{eq:crosseq-4}). Since $y_{(j)}$ is a WLS estimate of $Cx$, it can be written as $y_{(j)}=Cx+\upsilon_{(j)}$, where $\upsilon_{(j)}$ is the measurement noise. Then one has
\begin{eqnarray}
\label{eq:crosseq-9}
\upsilon_{(j)}=y_{(j)}-Cx
\end{eqnarray}
Substituting (\ref{eq:crosseq-5}) into (\ref{eq:crosseq-6}) and noting $y_{(j)}=\hat{z}_{(j)}$, one has by (\ref{eq:crosseq-9}) that
\begin{eqnarray}
\label{eq:crosseq-10}
\upsilon_{(j)}=[e^{\mathrm{T}}\bar{R}_{(j)}^{-1}e]^{-1}e^{\mathrm{T}}\bar{R}_{(j)}^{-1}\bar{\upsilon}_{(j)}
\end{eqnarray}
Since $\upsilon_{(j-1)}$ is uncorrelated with $\upsilon_{j+1}$, it follows from (\ref{eq:crosseq-10}) that $R_{(j)}=\mathrm{Cov}\{\upsilon_{(j)}\}=[R_{(j-1)}^{-1}+R_{j+1}^{-1}]^{-1}$, which is just the equation (\ref{eq:crosseq-3}). Moreover, it follows from the equation $R_{(j)}^{-1}=R_{(j-1)}^{-1}+R_{j+1}^{-1}$ that $R_{(j)}\leq R_{(j-1)}$ and $R_{(j)}\leq R_{j+1}$,
which leads to $R_{f}^{s}(k)\leq R_{i}(k)$, $\forall~~i\in\{1,\ldots,n\}$. The proof is thus completed.

When the fused measurement $y_{f}^{s}(k)$ is available, the estimator is able to produce an optimal state estimate by using $y_{f}^{s}(k)$ and applying a standard Kalman filter. An alternative approach to obtain the fused measurement is to apply the BMF method, which has been presented in \cite{Deng06}. In the BMF method, all the measurements in $Y(k)$ are fused all at a time and the fused measurement can be obtained by the WLS method, and it is given by
\begin{eqnarray}
\label{eq:crosseq-11}
R_{f}^{b}(k)&=&[\Sigma_{i=1}^{n}R_{i}^{-1}(k]^{-1}\\
\label{eq:crosseq-12}
y_{f}^{b}(k)&=&R_{f}^{b}(k)[\Sigma_{i=1}^{n}R_{i}^{-1}(k)y_{i}(k)]
\end{eqnarray}

{\bf Remark 1}. If the measurement equations in (\ref{eq:crosseq-2}) have different measurement matrices $C_{s,i}(k)$, $i\in\mathcal{N}_{s},~s\in \Phi$, and $C_{s,i}(k)$ can be decomposed as $C_{s,i}=M_{s,i}C$, $\forall~i\in\mathcal{N}_{s},~s\in \Phi$, where $C\in\Re^{q\times n_{x}}$, $M_{s,i}\in\Re^{q\times q}$ and $\sum\limits_{i\in\mathcal{N}_{s}} M_{s,i}^{\mathrm{T}}R_{i}^{-1}M_{s,i}$ is non-singular, then a similar SMF rule as given in Theorem 1 can be obtained by following some similar lines as in (\ref{eq:crosseq-5})-(\ref{eq:crosseq-10}).

It has been shown in \cite{Deng06} that the BMF is optimal in the sense that the noise variance of the fused measurement is minimal among all the fusion rules with matrix weights. The following theorem shows that the proposed SMF is equivalent to the BMF.

{\bf Theorem 2}. The SMF is equivalent to the BMF, i.e., $y_{f}^{s}(k)=y_{f}^{b}(k)$ and $R_{f}^{s}(k)=R_{f}^{b}(k)$.

{\bf Proof}. For $j=n-1$, one has by (\ref{eq:crosseq-3}) and (\ref{eq:crosseq-4}) that
\begin{eqnarray}
\label{eq:crosseq-13}
&&R_{f}^{s}=R_{(n-1)}=[R_{(n-2)}^{-1}+R_{n}^{-1}]^{-1}\\
\label{eq:crosseq-14}
&&y_{f}^{s}=y_{(n-1)}=R_{(n-1)}^{-1}[R_{(n-2)}^{-1}y_{(n-2)}+R_{n}^{-1}y_{n}]
\end{eqnarray}
Substituting the expressions of $R_{(n-2)}$ and $y_{(n-2)}$ into (\ref{eq:crosseq-13}) and (\ref{eq:crosseq-14}) yields
\begin{eqnarray}
\label{eq:crosseq-15}
R_{(n-1)}&=&[R_{(n-3)}^{-1}+R_{n-1}^{-1}+R_{n}^{-1}]^{-1}\\
\label{eq:crosseq-16}
y_{(n-1)}&=&R_{(n-1)}^{-1}[R_{(n-3)}^{-1}y_{(n-3)}+R_{n-1}^{-1}y_{n-1}\nonumber\\
&&~~~~~~~~~~~+R_{n}^{-1}y_{n}]
\end{eqnarray}
Following the similar procedures as in (\ref{eq:crosseq-15}) and (\ref{eq:crosseq-16}) for $j=n-3,n-4,\ldots,n$, one finally obtains
\begin{eqnarray}
\label{eq:crosseq-17}
&&R_{f}^{s}=\left[\Sigma_{i=1}^{n}R_{i}^{-1}\right]^{-1}=R_{f}^{b}\\
\label{eq:crosseq-18}
&&y_{f}^{s}=(R_{f}^{s})^{-1}\left[\Sigma_{i=1}^{n}R_{i}^{-1}y_{i}\right]=y_{f}^{b}
\end{eqnarray}
The proof is thus completed.

{\bf Remark 2}. Conventional approaches for the local estimation with multiple measurements either use the measurement augmentation (MA) or use the sequential Kalman (SK) estimation. It can be seen from Theorem 1 and equation (\ref{eq:crosseq-12}) that the fused measurement obtained by the SMF or BMF has the same dimension as each local measurement. Therefore, the SMF and BMF methods are more computationally efficient than the MA method. Moreover, in the SMF or BMF method, the estimate is obtained using the fused measurement and applying one step state update. However, in the SK estimation method, the estimate is obtained by sequentially applying a set of updates. Therefore, the SMF and BMF methods have lower computation complexity than the SK method. Specifically, Define the computational complexity of the estimation method as the number of multiplications and divisions
in the algorithm, and let $\delta_{sm}$, $\delta_{bm}$, $\delta_{ma}$, and $\delta_{sk}$ denote the
computational complexity of the SMF, BMF, MA and SK methods, respectively. Then, one has
$\delta_{sm}=(n_{x}^{2}+5n_{x}+8n-7)q^{2}+(4n_{x}^{2}+n_{x})q+2n_{x}^{3}+n_{x}^{2}$,
$\delta_{bm}=(n_{x}^{2}+5n_{x}+3n+3)q^{2}+(4n_{x}^{2}+n_{x})q+2n_{x}^{3}+n_{x}^{2}$,
$\delta_{ma}=(n_{x}^{2}+5n_{x}+8n-7)n^{2}q^{2}+(4n_{x}^{2}+n_{x})nq+2n_{x}^{3}+n_{x}^{2}$,
$\delta_{sk}=(n_{x}^{2}+5n_{x}+8n-7)nq^{2}+(4n_{x}^{2}+n_{x})nq+n(2n_{x}^{3}+n_{x}^{2})$. It can be seen that both $\delta_{sm}$ and $\delta_{bm}$ are of magnitude $\mathcal{O}(nq^{2})$, while $\delta_{ma}$ and $\delta_{sk}$ are of magnitudes $\mathcal{O}(n^{3}q^{2})$ and $\mathcal{O}(n^{2}q^{2})$, respectively.

%\begin{eqnarray}
%\label{eq:crosseq-xx}
%\delta_{sm}=(n_{x}^{2}+5n_{x}+8n-7)q^{2}+(4n_{x}^{2}+n_{x})q+2n_{x}^{3}+n_{x}^{2}\\
%\label{eq:crosseq-xx}
%\delta_{bm}=(n_{x}^{2}+5n_{x}+3n+3)q^{2}+(4n_{x}^{2}+n_{x})q+2n_{x}^{3}+n_{x}^{2}\\
%\label{eq:crosseq-xx}
%\delta_{ma}=(n_{x}^{2}+5n_{x}+8n-7)n^{2}q^{2}+(4n_{x}^{2}+n_{x})nq+2n_{x}^{3}+n_{x}^{2}\\
%\label{eq:crosseq-xx}
%\delta_{sk}=(n_{x}^{2}+5n_{x}+8n-7)nq^{2}+(4n_{x}^{2}+n_{x})nq+n(2n_{x}^{3}+n_{x}^{2})
%\end{eqnarray}

\section{Design of the SSF Estimators}
This section is devoted to the design of the state fusion estimator for each cluster. Suppose that $m$ local estimates $\hat{x}_{i}$, $i=1,2,\ldots,m$ are available for fusion at the cluster $\mathcal{N}_{s}$, $s\in \Phi$ over each estimation interval. To fuse the $m$ local estimates, one may apply some well developed batch state fusion (BSF) methods, such as the optimal fusion rule with matrix weights given in \cite{Sunsl04}. However, in the BSF method, all the local estimates are fused once at a time which usually involves computing the inverse of a high dimensional matrix, and may not be suitable for real-time applications. In what follows, a SSF method with matrix weights will be developed.

Without loss of generality, suppose that the $m$ local estimates arrive at the head of cluster $\mathcal{N}_{s}$ in time order as $\hat{x}_{1}$, $\hat{x}_{2}$, $\ldots$, $\hat{x}_{m}$. In the SSF method, the cluster head fuses the local estimates one by one according to the time order. Denote the $j$th fused estimate as $\hat{x}_{(j)}$, then $j\in\{1,2,\ldots,m-1\}$, $\hat{x}_{(0)}=\hat{x}_{1}$ and $\hat{x}_{(j)}=f_{s}(\hat{x}_{(j-1)},\hat{x}_{j+1})$, where $f_{s}$ is the state fusion rule to be designed. The SSF with matrix weights are presented in the following theorem.

{\bf Theorem 3}. Let $\hat{x}_{i}$, $i=1,2,\ldots,m$ be unbiased estimates of the state of system (\ref{eq:crosseq-1}) and $P_{i}$ be the estimation error variance matrix of $\hat{x}_{i}$. Then, the SSF estimator with matrix weights is given by
\begin{eqnarray}
\label{eq:crosseq-19}
\hat{x}_{(j)}&=&\Delta_{1,(j)}\hat{x}_{(j-1)}+\Delta_{2,(j)}\hat{x}_{j+1}\\
\label{eq:crosseq-20}
P_{(j)}&=&[e^{\mathrm{T}}\Omega_{(j)}^{-1}e]^{-1}
\end{eqnarray}
where $j=1,2,\ldots,m-1$, $P_{(j)}$ is the estimation error variance matrix of $\hat{x}_{(j)}$, $\hat{x}_{(0)}=\hat{x}_{1}$, $P_{(0)}=P_{1}$ and the optimal matrix weights $\Delta_{1,(j)}$ and $\Delta_{2,(j)}$ are computed as
\begin{eqnarray}
\label{eq:crosseq-21}
&&\left[\begin{array}{c}\Delta_{1,(j)}\\\Delta_{2,(j)}\end{array}\right]=\Omega_{(j)}^{-1}e[e^{\mathrm{T}}\Omega_{(j)}^{-1}e]^{-1}\\
\label{eq:crosseq-22}
&&\Omega_{(j)}=\left[\begin{array}{cc}P_{(j-1)}&P_{(j-1),j+1}\\\ast&P_{j+1}\end{array}\right]\\
\label{eq:crosseq-23}
&&e=[I~~I]^{\mathrm{T}}
\end{eqnarray}
where $P_{(j-1),j+1}=\mathbf{E}\{\tilde{x}_{(j-1)}\tilde{x}_{j+1}^{\mathrm{T}}\}$ is the cross covariance matrix of $\hat{x}_{(j-1)}$ and $\hat{x}_{j+1}$, and is computed as
\begin{eqnarray}
\label{eq:crosseq-24}
&&P_{(j-1),j+1}(k)=\sum_{d=1}^{j}\prod_{l=1}^{j-d+1}\Delta_{1,(j-l)}\Delta_{2,(d-1)}P_{d,j+1}(k)\\
\label{eq:crosseq-25}
&&P_{j,d}(k)=[I-K_{j}(k)C(k)][A(k-1)P_{j,d}(k-1)\nonumber\\
&&~~~~~~~~~~\times A^{\mathrm{T}}(k-1)+B(k-1)Q_{\omega}B^{\mathrm{T}}(k-1)]\nonumber\\
&&~~~~~~~~~~\times [I-K_{d}(k)C(k)]^{\mathrm{T}}
\end{eqnarray}
where $\Delta_{1,(0)}=\Delta_{2,(0)}=I$ and $P_{(0),2}=P_{1,2}$. The fused state estimate and its error variance are finally given by $\hat{x}_{f}^{s}=\hat{x}_{(m-1)}$ and $P_{f}^{s}=P_{(m-1)}$, and one has $P_{f}^{s}\leq P_{i}$, $i=1,2,\ldots,m$, that is, the precision of the SSF estimator is higher than each local estimator.

{\bf Proof}. According to the SSF rule, $\hat{x}_{(j-1)}$ and $\hat{x}_{j+1}$ are fused in the $j$th fusion. Then, by applying the fusion rule with matrix weights
as presented in Theorem 1 of \cite{Sunsl04}, the optimal fusion in the linear minimum variance sense is given by (\ref{eq:crosseq-19}) and (\ref{eq:crosseq-20}), where $\Delta_{1,(j)}+\Delta_{2,(j)}=I$. In what follows, it will be shown that the cross-covariance $P_{(j-1),j+1}$ satisfies the equation (\ref{eq:crosseq-24}). Note that for $j=1$, one has $P_{(0),2}=P_{1,2}$. For $j\geq2$ and $t\in\{j-1,j-2,\ldots,1\}$ one has
\begin{eqnarray}
\label{eq:crosseq-26}
&&\tilde{x}_{(t)}=x-\hat{x}_{(t)}\\
\label{eq:crosseq-27}
&&\hat{x}_{(t)}=\Delta_{1,(t)}\hat{x}_{(t-1)}+\Delta_{2,(t)}\hat{x}_{t+1}
\end{eqnarray}
Substituting (\ref{eq:crosseq-27}) into (\ref{eq:crosseq-26}) and taking the relation $\Delta_{1,(t)}+\Delta_{2,(t)}=I$ into consideration yields
\begin{eqnarray}
\label{eq:crosseq-28}
\tilde{x}_{(t)}&=&(\Delta_{1,(t)}+\Delta_{2,(t)})x-\Delta_{1,(t)}\hat{x}_{(t-1)}-\Delta_{2,(t)}\hat{x}_{t+1}\nonumber\\
&=&\Delta_{1,(t)}\tilde{x}_{(t-1)}+\Delta_{2,(t)}\tilde{x}_{t+1}
\end{eqnarray}
It follows from (\ref{eq:crosseq-28}) that
\begin{eqnarray}
\label{eq:crosseq-29}
P_{(t),j+1}&=&\mathbf{E}\{\tilde{x}_{(t)}\tilde{x}_{j+1}^{\mathrm{T}}\}\nonumber\\
&=&\Delta_{1,(t)}P_{(t-1),j+1}+\Delta_{2,(t)}P_{t+1,j+1}
\end{eqnarray}
Applying (\ref{eq:crosseq-29}) recursively for $t=j-1,j-2,\ldots,1$ leads to equation (\ref{eq:crosseq-24}).
Denote the fused measurement in the $j$th cluster as $y_{f,j}^{s}$, then $y_{f,j}^{s}$ can be written as
\begin{eqnarray}
\label{eq:crosseq-30}
y_{f,j}^{s}=Cx+\upsilon_{f,j}^{s},~j\in Z_{s}
\end{eqnarray}
where $\upsilon_{f,j}^{s}$ is the fused measurement noise. Then, by (\ref{eq:crosseq-1}), (\ref{eq:crosseq-30}) and the standard Kalman filter, one has
\begin{eqnarray}
\label{eq:crosseq-31}
\tilde{x}_{j}(k)&=&x(k)-\hat{x}_{j}(k)\nonumber\\
&=&[I-K_{j}(k)C(k)][A(k-1)\tilde{x}_{j}(k-1)\nonumber\\
&&+B(k-1)\omega(k-1)]
-K_{j}(k)\upsilon_{f,j}^{s}(k)
\end{eqnarray}
By definition, the cross covariance is given by
\begin{eqnarray}
\label{eq:crosseq-32}
P_{j,d}(k)=\mathbf{E}\{\tilde{x}_{j}(k)\tilde{x}_{d}^{\mathrm{T}}(k)\}, j\neq d,~j,d\in Z_{s}
\end{eqnarray}
Substituting (\ref{eq:crosseq-31}) into (\ref{eq:crosseq-32}) and taking into account the relations $\tilde{x}_{j}(k-1)\perp \omega(k-1)$, $\tilde{x}_{j}(k-1)\perp \upsilon_{f,j}^{s}(k)$, $\omega(k-1)\perp\upsilon_{f,j}^{s}(k)$ and $\upsilon_{f,j}^{s}(k)\perp\upsilon_{f,d}^{s}(k)$, one obtains equation (\ref{eq:crosseq-25}).

Moreover, by the fusion rule given in \cite{Sunsl04}, one has
\begin{eqnarray}
\label{eq:crosseq-33}
P_{(j)}\leq P_{(j-1)},~P_{(j)}\leq P_{j+1}
\end{eqnarray}
By applying (\ref{eq:crosseq-33}) recursively for $j=1,2,\ldots,m-1$, one obtains $P_{f}^{s}=P_{(m-1)}\leq P_{i}$, $i\in Z_{s}$. The proof is thus completed.

{\bf Remark 3}. It can be seen from (\ref{eq:crosseq-22}) that the matrix $\Omega_{(j)}$ has dimension $2n_{x}\times2n_{x}$, while it has dimension $mn_{x}\times mn_{x}$ in the BSF. Therefore, the SSF is much more computational efficient than BSF. Specifically, let $\delta_{bs}$ and $\delta_{ss}$ denote the
computational complexity of the BSF and SSF methods, respectively. Then, one has $\delta_{bs}=5n_{x}^{2}m^{2}+(n_{x}^{3}+n_{x}^{2})m$ and $\delta_{ss}=(2n_{x}^{3}+22n_{x}^{2})m-2n_{x}^{3}-22n_{x}^{2}$. Therefore, the magnitudes of $\delta_{bs}$ and $\delta_{ss}$ are $\mathcal{O}(m^{2})$ and $\mathcal{O}(m)$, respectively.

The fused estimate $\hat{x}_{f}^{s}$ given by Theorem 3 has the following property.

{\bf Theorem 4}. The fused estimate $\hat{x}_{f}^{s}$ is an unbiased estimate of the system state, and satisfies
\begin{eqnarray}
\label{eq:crosseq-34}
\left\{\begin{array}{ccc}\hat{x}_{f}^{s}&=&\Delta\hat{x}\\\Delta I_{o}&=&I\end{array}\right.
\end{eqnarray}
where $\hat{x}=[\hat{x}_{1}^{\mathrm{T}}~\cdots~\hat{x}_{m}^{\mathrm{T}}]^{\mathrm{T}}$,
$I_{o}=[\underbrace{I~\cdots~I}\limits_{m}]^{\mathrm{T}}$ and
\begin{eqnarray*}
\Delta=\left[\prod_{l=1}^{m-1}\Delta_{1,(m-l)}~~\prod_{l=1}^{m-2}\Delta_{1,(m-l)}\Delta_{2,(1)}~~\cdots\right.\\
\left.\Delta_{1,(m-1)}\Delta_{2,(m-2)}~~\Delta_{2,(m-1)}\right]
\end{eqnarray*}

{\bf Proof}. Note that (\ref{eq:crosseq-19}) is a recursive equation for computing $\hat{x}_{(j)}$ with respect to the variable $j$. Therefore, substituting the expression of $\hat{x}_{(j-1)}$ into that of $\hat{x}_{(j)}$ for $j=1,2,\ldots,m-1$ yields $\hat{x}_{f}^{s}=\hat{x}_{(m-1)}=\Delta\hat{x}$. By the definition of $I_{o}$, one has
\begin{eqnarray}
\label{eq:crosseq-35}
\Delta I_{o}=\sum_{j=1}^{m}\prod\limits_{l=1}^{m-j}\Delta_{1,(m-l)}\Delta_{2,(j-1)}
\end{eqnarray}
Denote $D_{m-j+1}=\prod\limits_{l=1}^{m-j}\Delta_{1,(m-l)}\Delta_{2,(j-1)}$, $j=1,\ldots,m$. Then, one has by (\ref{eq:crosseq-35}) that
\begin{eqnarray}
\label{eq:crosseq-36}
\Delta I_{o}=\sum_{j=1}^{m}D_{j}
\end{eqnarray}
Since $\Delta_{1,(j)}+\Delta_{2,(j)}=I$, $j=1,2,\ldots,m-1$, one has
\begin{eqnarray}
\label{eq:crosseq-37}
D_{1}+D_{2}&=&\Delta_{1,(m-1)}\Delta_{2,(m-2)}+\Delta_{2,(m-1)}\nonumber\\
&=&\Delta_{1,(m-1)}\Delta_{2,(m-2)}+I-\Delta_{1,(m-1)}\nonumber\\
&=&I+\Delta_{1,(m-1)}[\Delta_{2,(m-2)}-I]\nonumber\\
&=&I-\Delta_{1,(m-1)}\Delta_{1,(m-2)}
\end{eqnarray}
Then, it follows from (\ref{eq:crosseq-37}) that
\begin{eqnarray}
\label{eq:crosseq-38}
&&D_{1}+D_{2}+D_{3}\nonumber\\
&&=I-\Delta_{1,(m-1)}\Delta_{1,(m-2)}+\Delta_{1,(m-1)}\Delta_{1,(m-2)}\Delta_{2,(m-3)}\nonumber\\
&&=I+\Delta_{1,(m-1)}\Delta_{1,(m-2)}[\Delta_{2,(m-3)}-I]\nonumber\\
&&=I-\Delta_{1,(m-1)}\Delta_{1,(m-2)}\Delta_{1,(m-3)}
\end{eqnarray}
Following the similar lines as in (\ref{eq:crosseq-37}) and (\ref{eq:crosseq-38}), one obtains
\begin{eqnarray}
\label{eq:crosseq-39}
\sum_{j=1}^{m-1}D_{j}=I-\prod_{l=1}^{m-1}\Delta_{1,(m-l)}
\end{eqnarray}
Since $\Delta_{2,(0)}=I$, one has
\begin{eqnarray}
\label{eq:crosseq-40}
D_{m}=\prod_{l=1}^{m-1}\Delta_{1,(m-l)}\Delta_{2,(0)}=\prod_{l=1}^{m-1}\Delta_{1,(m-l)}
\end{eqnarray}
Then, it follows from (\ref{eq:crosseq-36}), (\ref{eq:crosseq-39}) and (\ref{eq:crosseq-40}) that
\begin{eqnarray}
\label{eq:crosseq-41}
\Delta I_{o}=\sum_{j=1}^{m}D_{j}=\sum_{j=1}^{m-1}D_{j}+D_{m}=I
\end{eqnarray}
Since $\hat{x}_{i}$, $i=1,2,\ldots,m$ are unbiased estimates of $x$, one has by the fact $\Delta I_{o}=I$ that $\mathbf{E}\{x-\hat{x}_{f}^{s}\}=\mathbf{E}\{\Delta I_{o}x-\Delta\hat{x}\}=\Delta\mathbf{E}\{I_{o}x-\hat{x}\}=0$, that is, $\mathbf{E}\{x\}=\mathbf{E}\{\hat{x}_{f}^{s}\}$. Therefore, $\hat{x}_{f}^{s}$ is an unbiased estimate of the state $x$. The proof is thus completed.

{\bf Remark 4}. It can be seen from Theorem 4 that the proposed estimate $\hat{x}_{f}^{s}$ is also a linear combination of all the local estimates with matrix weights. Thus, the estimator given in Theorem 4 can be regarded as another form of BSF estimator. However, the weighting matrices in the BSF estimator given in Theorem 4 have much lower dimensions than those in the BSF estimator as presented in \cite{Sunsl04}.

\section{Simulations}
Consider a maneuvering target tracking system, where the
target moves in one direction, and its position and velocity
evolve according to the state-space model (\ref{eq:crosseq-1}) with
\begin{eqnarray}
\label{eq:crosseq-42}
A(k)=\left[\begin{array}{cc}1&h\\0&1\end{array}\right],
B(k)=\left[\begin{array}{c}h^{2}/2\\h\end{array}\right]
\end{eqnarray}
where $h$ is the sampling period. The state is
$x(k)=[x_{p}^{\mathrm{T}}(k)~~x_{v}^{\mathrm{T}}(k)]^{\mathrm{T}}$, where
$x_{p}(k)$ and $x_{v}(k)$ are the position and velocity of
the maneuvering target at time $k$, respectively. The variance of the process noise $\omega(k)$ is $q=1$. The position is measured by the sensors, and the measurement matrix is $C=[1~0]$.
The initial state of system (\ref{eq:crosseq-1}) is $x_{0}=[1~~0.5]^{\mathrm{T}}$, and we take $h=0.5s$ in the simulation. A group of sensors are deployed to monitor the target, and the sensor network is divided into three clusters, namely, $\mathcal{N}_{1}$ with 10 sensors, $\mathcal{N}_{2}$ with 8 sensors, and $\mathcal{N}_{3}$ with 6 sensors. There is a CH in each cluster, and the CH collects measurements from its cluster to generate a local estimate of the system state. Monte Carlo simulations will be carried out and the root mean-square error
$\mathrm{RMSE}=\sqrt{\frac{1}{L}\sum_{i=1}^{L}(x_{p}^{i}-\hat{x}_{p}^{i})^{2}}$
is used to evaluate the estimation performance of the estimators, where $L=1000$ is the number of Monte Carlo simulation runs. The estimate of the initial state is set as $\hat{x}_{0}=[2~~1]^{\mathrm{T}}$.

\begin{figure}
\begin{center}
\includegraphics[height=6cm]{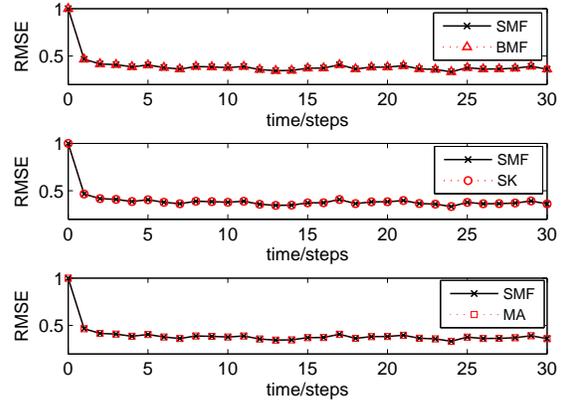}
\caption{Equivalence of the BMF estimator and the SMF, SK and MA estimator}
\label{fig3}
\end{center}
\end{figure}

\begin{figure}
\begin{center}
\includegraphics[height=4cm]{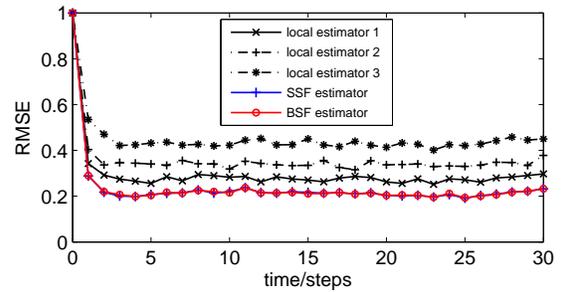}
\caption{Equivalence of the BSF estimator and the SSF estimator}
\label{fig4}
\end{center}
\end{figure}

The local estimates in cluster $\mathcal{N}_{1}$ using the proposed SMF method, the BMF, SK and MA methods are shown in Fig. \ref{fig3}. It can be seen from Fig. \ref{fig3} that the SMF, BMF, SK and MA methods provide the same estimation precision. Now, consider the state fusion in cluster $\mathcal{N}_{1}$. In the state fusion stage, the CH in $\mathcal{N}_{1}$ collects local estimates from the other two CHs and generates the fused state estimates using the proposed BSF and SSF methods. It can be seen from Fig. \ref{fig4} that the estimation performance is improved by fusing the local estimates, and the proposed SSF method is equivalent to the BSF method in achieving the same estimation precision.

\section{Conclusions}
Some sequential fusion estimators have been developed in this paper for distributed estimation in clustered sensor networks. It is shown that the sequential fusion methods have the same estimation performance as the batch fusion one but have lower computational complexity than conventional approaches, such as the batch fusion estimation, sequential Kalman estimation and that based on measurement augmentation. Therefore, the proposed methods are more appropriate for real-time applications and are convenient to handle asynchronous and delayed information.

%We are planning to expand our research to the case of multi-sensor fusion estimation for nonlinear systems in the future.

%\begin{ack}                               % Place acknowledgements
%Partially supported by xxx.  % here.
%\end{ack}

\bibliographystyle{plain}        % Include this if you use bibtex
\bibliography{autosam}           % and a bib file to produce the

\begin{thebibliography}{99}     % Otherwise use the

\bibitem{IlicMD10}
Ilic, M.D., Xie, L., Khan, U.A., Moura, J.M.F. (2010). Modeling of future cyber-physical energy systems for distributed
sensing and control. {\em IEEE Transactions on Systems Man and Cybernetics, Part A--Systems and Humans}, {\em 40}(4), 825--838.

\bibitem{CaoXH14}
Cao, X.H., Cheng, P., Chen, J.M., Ge, S.Z., Cheng, Y., Sun Y.X. (2014). Cognitive radio based state estimation in cyber-physical systems. {\em IEEE Journal on Selected Areas in Communications}, {\em 32}(3): 489--502.

\bibitem{ChenJM13}
Chen, J.M., Li J.K., Lai, T.H. (2013). Trapping mobile targets in wireless sensor networks: an energy-efficient perspective. {\em IEEE Transactions on Vehicular Technology}, {\em 62}(7): 3287--3300.

\bibitem{HeX11}
He, X., Wang, Z.D., Ji, Y.D., Zhou, D.H. (2011). Robust fault detection for networked systems with distributed sensors. {\em IEEE Transactions on Aerospace and Electronic Systems}, {\em 47}(1): 166--177.

\bibitem{OkaA10}
Oka, A., Lampe, L. (2010). Distributed target tracking using signal strength measurements by a wireless sensor network. {\em IEEE Journal on Selected Areas in Communications}, {\em 28}(7): 1006--1015.

\bibitem{DongHL12}
Dong, H.L., Wang, Z.D., Gao, H.J. (2012). Distributed filtering for a class of time-varying systems over sensor networks with quantization errors and successive packet dropouts. {\em IEEE Transactions on Signal Processing}, {\em 60}(6): 3164--3173.

\bibitem{MillanP13}
Millan, P., Orihuela, L., Vivas, C., Rubio, F.R., Dimarogonas, D.V., Johansson, K.H. (2013). Sensor-network-based robust distributed control and estimation. {\em Control Engineering Practice}, {\em 21}(9): 1238--1249.

\bibitem{SongHY14}
Song, H.Y., Zhang, W.A., Yu, L. (2014). Hierarchical fusion in clustered sensor networks with asynchronous local estimates. {\em IEEE Signal Processing Letters}, {\em 21}(12): 1506--1510.

\bibitem{Roeckerja88}
Roecker, J.A., McGillem, C.D. (1988) Comparison of two-sensor tracking methods based on state vector fusion and measurement fusion. {\em IEEE Transactions on Aerospace and Electronic Systems}, {\em 24}(4), 447--449.

\bibitem{BarShalom95}
Bar-Shalom, Y., and Li, X. R. (1995). Multitarget-multisensor tracking:
Principles and techniques. Storrs, CT: YBS Publishing.

\bibitem{Sunsl04}
Sun, S.L., Deng, Z.L. (2004). Multi-sensor optimal information fusion
Kalman filter. {\em Automatica}, {\em 40}(6), 1017--1023.

\bibitem{Deng06}
Deng, Z.L. (2006). On functional equivalence of two measurement fusion methods. {\em
Control Theory and Applications}, {\em 23}(2), 319--323.

\bibitem{Songeb07}
Song, E.B., Zhu, Y.M., Zhou, J., You, Z.S. (2007). Optimal Kalman filtering fusion with cross-correlated sensor noises.
{\em Automatica}, {\em 43}(8), 1450--1456.

\bibitem{Julier09}
Julier, S.J., Uhlman, J.K. (2009). General decentralized data fusion with covariance intersection, in: M.E. Liggins, D.L. Hall, J. Llinas (Eds.), Handbook of
multisensor data fusion, Second ed., Theory and Practice, CRC Press.

\bibitem{YYhu10}
Hu, Y.Y., Duan, Z.S., and Zhou, D.H. (2010). Estimation fusion with general asynchronous multi-rate sensors.
{\em IEEE Trans. Aerosp. Electron. Syst.}, {\em 46}(4), 2090--2102.

\bibitem{Zhang14}
Zhang, W.A., Liu, S., and Yu, L. (2014). Fusion estimation for sensor networks with uniform estimation rates. {\em IEEE Trans. Circuits Syst.--I: Reg. Papers}, {\em 61}(5), 1485--1498.

\bibitem{XiaYQ09}
Xia, Y.Q., Shang, J.H., Chen, J., and Liu, G.P. (2009). Networked data fusion with packet losses and variable delays. {\em IEEE Transactions on Systems, Man and Cybernetics--Part B}, {\em 39}(5), 1107--1120.

\bibitem{XingZR16}
Xing, Z.R., Xia, Y.Q. (2016). Distributed federated Kalman filter fusion over multi-sensor unreliable networked systems. {\em IEEE Transactions on Circuits and Systems I--Regular Papers}, {\em 63}(10), 1714--1725.

\bibitem{YanLP13}
Yan, L.P., Li, X.R., Xia, Y.Q., and Fu, M.Y. (2013). Optimal sequential and distributed fusion for state estimation in
cross-correlated noise. {\em Automatica}, {\em 49}, 3607--3612.

\bibitem{Deng12}
Deng, Z.L., Zhang, P., Qi, W.J., Liu, J.F., and Gao, Y. (2012). Sequential covariance intersection fusion Kalman filter. {\em Information Sciences}, {\em 189}, 293--309.

\end{thebibliography}
                                 % bibliography (preferred). The
                                 % correct style is generated by
                                 % Elsevier at the time of printing.

%\appendix

\end{document}